# ECOLOGICAL EQUILIBRIUM FOR RESTRAINED BRANCHING RANDOM WALKS

By Daniela Bertacchi, Gustavo Posta and Fabio Zucca

*Università di Milano–Bicocca, Politecnico di Milano and Politecnico di Milano*

We study a generalized branching random walk where particles breed at a rate which depends on the number of neighboring particles. Under general assumptions on the breeding rates we prove the existence of a phase where the population survives without exploding. We construct a nontrivial invariant measure for this case.

**1. Introduction.** Scientists have been studying models for the evolution of a population since the end of the 19th century, starting from the branching process introduced by Galton and Watson in 1875 [5]. The need for more realistic models has led to the introduction of a spatial structure: the branching random walk and the contact process (briefly, BRW and CP, resp.) are perhaps the most natural generalizations. In the BRW model each individual has a fixed position on a connected graph, for example, the integer lattice $\mathbb{Z}^d$, and an exponential lifespan of parameter 1 during which it breeds on neighboring sites according to a Poisson process of intensity $\lambda > 0$. The number of individuals allowed per site is unbounded. Requiring that a site can be occupied by at most one individual, one obtains the CP. Both these processes exhibit two possible behaviors: starting from a finite population, either the population faces almost sure extinction (subcritical behavior), or it survives with a positive probability (supercritical behavior). In the supercritical case the BRW's population grows indefinitely and the mean density of the population diverges. For the contact process, obviously there is no divergence of the mean density of the population because this quantity is a priori bounded. In the supercritical phase, the CP has two invariant extremal measures (see [8]). It is known that there exists a critical value of $\lambda$ separating the two behaviors: if $\lambda$ is smaller than the critical parameter, the process exhibits the subcritical behavior, while, for larger $\lambda$, it exhibits the supercritical one. We denote by $\lambda_{\text{BRW}}$ and $\lambda_{\text{CP}}$ the critical parameters of the BRW and of the CP respectively.

The observation of natural environments suggests to remove any a priori bound on the number of individuals allowed per site and to introduce a self-regulating mechanism on the birth rates, which should provide a surviving though nonexploding population. Indeed, some ecological systems seem to be in a sort of equilibrium where the density of a population neither tends to zero nor to infinity. One may argue that we could be observing a subcritical or supercritical system during a too short time span, nevertheless, it seems natural to try to translate into mathematical terms the competition for resources (see, e.g., the discussion in [7]). Other authors have introduced models for self-regulating populations. For instance, in the case of a population living on a continuous and homogeneous space, Bolker and Pacala [2] studied a process where the death rates depend on the local density centered on the father. A slightly different model was considered in [4] where the reproduction rate depends on the local density centered on the father. The main technical tools are moment equations and stochastic differential equations respectively. A different approach is carried out in [3] where the population has no spatial structure and each individual can be affected by a gene mutation at birth; the evolution is studied as a Markov process in the trait space.

We introduce a self-regulating mechanism where the birth rate is a decreasing function of the local density at the location where the offspring would live. Moreover, noting that the spatial structure of the interaction between individuals in a biological population might be irregular, we study a population on a discrete (possibly nonhomogeneous) space. To this aim, we consider the following model, which we call *restrained branching random walk* (RBRW briefly). Consider an infinite connected graph $X$ with bounded geometry (i.e., the number of neighbors of the vertices is bounded, e.g., $\mathbb{Z}^d$) as the environment where the population lives and let $\eta(x)$ be the number of individuals living at the site $x \in X$. The lifespan of each individual is an exponential random variable of mean 1. During its lifetime each individual tries to reproduce following a Poisson process of intensity $\lambda$. Every time the clock associated to the Poisson process rings, the individual tries to send an offspring to a randomly chosen target neighboring site. The target neighboring site is chosen using the transition matrix $P = (p(x,y))_{x,y \in X}$ of a nearest neighbor random walk on $X$, for example, the simple random walk on $\mathbb{Z}^d$. Call the target site $y$. The reproduction on $y$ is effective only with probability $c(\eta(y))/\lambda$, where $c : \mathbb{N} \to \mathbb{R}^+$ is a nonincreasing and nonnegative function with $c(0) = \lambda$. In this case the population living at $y$ increases by one individual, otherwise nothing happens.

Observe that the process described above is a Markov process and includes the BRW and the CP as special cases ($c \equiv \lambda$ and $c = \lambda \mathbb{1}_{\{0\}}$, resp.). The formal construction of this process is carried out in Section 3, where the existence of a Markov process $\{\eta_t\}_{t \geq 0}$ with state space $\Omega \subset \mathbb{N}^X$ is proven. In general $\Omega$ is smaller than $\mathbb{N}^X$ because we can only consider configurations $\eta$ such that $\eta(x)$ does not diverge too fast when $x$ goes to infinity (see Section 2 for more details). We prove that $\{\eta_t\}_{t \geq 0}$ has different behaviors depending on $c(0)$, $c(+\infty) := \lim_{k \to +\infty} c(k)$ and on the transition kernel $p(x, y)$ (see Proposition 4.1 for the complete statement).

PROPOSITION 1.1. *Let $P^n = (p^{(n)}(x, y))_{x,y \in X}$ be the nth power of the transition matrix $P$ and $\{\eta_t\}_{t \geq 0}$ be the RBRW described above. Let us define $\rho := \limsup_{n \to \infty} \sqrt[n]{p^{(n)}(x, y)}$ and $\theta := \lim_{n \to \infty} \sqrt[n]{\sup_x \sum_y p^{(n)}(y, x)}$ (notice that $\theta \geq \rho$):*

(i) *If $c(0) < 1/\rho$, then $\lim_{t \to +\infty} \mathbb{E}^\eta[\eta_t(x)] = 0$ for any finite $\eta \in \Omega$, $x \in X$;*
(ii) *If $c(0) > \lambda_{\mathrm{CP}}$, then $\lim_{t \to +\infty} \mathbb{E}^\eta[\eta_t(x)] > 0$ for any $\eta \in \Omega \setminus \{\mathbf{0}\}$, $x \in X$ and $\mathbb{P}^\eta(\limsup_{t \to \infty} \eta_t(x) > 0) > 0$;*
(iii) *If $c(+\infty) > 1/\rho$, then $\lim_{t \to +\infty} \mathbb{E}^\eta[\eta_t(x)] = +\infty$ for any $\eta \in \Omega \setminus \{\mathbf{0}\}$;*
(iv) *If $c(+\infty) < 1/\theta$, then $\limsup_{t \to +\infty} \mathbb{E}^\eta[\eta_t(x)] < +\infty$ uniformly for any bounded $\eta \in \Omega$, $x \in X$.*

The critical parameters $\rho$ and $\theta$ (and $\lambda_{\mathrm{CP}}$ as well) depend only on $P$ (hence, in the case of the simple random walk, on the geometry of the graph); in particular, $\lambda_{\mathrm{BRW}} = 1/\rho$ (see [15] and [1]). We discuss further details in Section 2.1. The proof of Proposition 1.1 is quite simple and essentially based on coupling techniques with the CP and the BRW (and on explicit estimates on the moments of the BRW with immigration; see Lemma 3.3).

Notice that, given a bounded initial state $\eta \in \Omega \setminus \{\mathbf{0}\}$, if $c(0)$ is sufficiently large and $c(+\infty)$ is sufficiently small, then by (ii) of Proposition 1.1, the population has a positive probability to survive indefinitely, while by (iv), almost surely, it does not explode. This is the ecological equilibrium phase we are looking for. It is quite natural to wonder if there is a stationary distribution for the population in this case. We prove that this is the case (see Theorem 4.3 for the complete statement).

THEOREM 1.2. *Let $\{\eta_t\}_{t \geq 0}$ be the RBRW described above and assume that $c$ is such that $c(0) > \lambda_{\mathrm{CP}}$ and $c(+\infty) < 1/\theta$. Then there exists a nontrivial probability measure $\mu$ on $(\Omega, \mathcal{B}(\Omega))$ which is invariant for $\{\eta_t\}_{t \geq 0}$.*

We construct this invariant measure as a limit of invariant measures of processes where the number of individuals per site is bounded.

Important examples are the RBRW on the $d$-dimensional lattice and on the homogeneous tree of degree $n+1$ (both endowed with the transition matrix of the simple random walk).

If $X = \mathbb{Z}^d$ and $P$ is the simple random walk, we have that $\rho = \theta = 1$. So (i) of Proposition 1.1 implies that the population dies out when $c(0)$ is smaller than the death rate 1, while (iii) states that the mean density of the population explodes when $c(+\infty)$ is larger than 1. Moreover, the system can reach the ecological equilibrium if $c(+\infty) < 1$ and $c(0) > \lambda_{\mathrm{CP}}$ (i.e., the critical parameter of the CP on $\mathbb{Z}^d$). In this case the stationary measure $\mu$ given by Theorem 1.2 is translation invariant.

If $X = \mathbb{T}_{n+1}$ (the homogeneous tree where the degree of each vertex is $n+1$) and $P$ is the simple random walk, we have $\rho = 2\sqrt{n}/(n+1) < 1 = \theta$, whence, to ensure ecological equilibrium we require $c(+\infty) < 1$ and $c(0) > \lambda_{\mathrm{CP}}$ (i.e., the critical parameter of the CP on $\mathbb{T}_{n+1}$), while, for the almost sure extinction, it is sufficient that $c(0) < (n+1)/2\sqrt{n}$. The stationary measure $\mu$ given by Theorem 1.2 is translation invariant in this case as well. Considering different random walks on $\mathbb{T}_{n+1}$ leads to different values for $\theta$ (see Example 5.1).

One may wonder how the two parameters $\rho$ and $\theta$ come to surface: the analysis of the two examples above, shows that on general graphs the behavior of interacting particle systems can be different than on $\mathbb{Z}^d$. For instance, it is known that on some fast growing graphs there is the so-called *weak phase* (see, e.g., [11] and [9], Part I, Chapter 4 for the CP on trees, [15] and [1] for the BRW on graphs and [12] for the BRW on Galton–Watson trees): the population can survive by drifting to infinity and leaving eventually any site.

On $\mathbb{Z}^d$, $\lambda_{\mathrm{BRW}} = 1$ and, in the subcritical phase, starting from a bounded $\eta$ with at least one individual per site, with probability one, there is no extinction and the expected number of individuals at a fixed site is a bounded function of the time $t$.

On a general graph the subcritical phase of the BRW is further subdivided: $\lambda_{\mathrm{BRW}} = 1/\rho$, but only if $\lambda < 1/\theta$ one can ensure that the expected number of individuals at a fixed site is a bounded function of the time $t$ (starting from a bounded $\eta$ with at least one individual per site). Indeed, if $\lambda \in (1/\theta, 1/\rho)$ and $\eta$ is as above, then there are examples where the expected number of individuals at a fixed site diverges as $t$ goes to infinity (e.g., if $P$ is the simple random walk on a homogeneous tree: see Example 5.1).

We give here a brief outline of the paper. In Section 2 we give the definitions needed in the sequel and we introduce the generator of the process. The construction of the RBRW is carried out in Section 3. Since the state space of this process is not locally compact, the classical approach of Hille–Yosida cannot be used: we follow the ideas of [10]. Some of the results we prove are obtained via a coupling argument (see Proposition 3.5) with particular

BRWs (with immortal particles). Furthermore, we give explicit estimates of some moments of these processes (see Lemma 3.3). In Section 4 we prove our main results: Proposition 1.1 and Theorem 1.2 (see Proposition 4.1 and Theorem 4.3, resp.). Section 5 is devoted to final remarks, examples and open questions.

## 2. Preliminaries.

2.1. *Graph geometry and random walks.* Let $X$ be a connected, nonoriented graph, with bounded geometry (i.e., the number of neighbors of a vertex $x$, called degree of $x$, is uniformly bounded on $X$); denote by $D$ the maximum degree of vertices on $X$. Let $P = (p(x,y))_{x,y \in X}$ be a stochastic matrix (although one can apply easily our methods to a substochastic $P$) such that $p(x,y) > 0$ if and only if $x$ and $y$ are neighbors (we write $x \sim y$ in this case). For any $\Lambda \subset X$, let $\Lambda^\circ := \{x \in \Lambda : \forall y \sim x, y \in \Lambda\}$ be the *interior* of $\Lambda$ and let

$$p_\Lambda(x,y) := \begin{cases} p(x,y), & \text{if } x,y \in \Lambda, \\ 0, & \text{otherwise.} \end{cases}$$

The two parameters $\rho$ and $\theta$, associated to $P$, play a crucial role in distinguishing between different behaviors of the RBRW (see Proposition 1.1). Recall that in Section 1 we defined $\rho := \limsup_{n \to \infty} \sqrt[n]{p^{(n)}(x,y)}$: this is usually called the *convergence parameter* and it is independent of $x, y \in X$ (see [14]).

Consider the space of (infinite) matrices endowed with the norm $\|A\| := \sup_x \sum_y |a_{xy}|$. As usual, each matrix with a finite norm can be identified with a linear continuous operator from $l^\infty(X)$ into itself. Let $\theta(A)$ be the spectral radius of the operator $A$; note that $\theta(A) = \lim_{n \to \infty} \|A^n\|^{1/n}$ (see, e.g., [13], Theorem 18.9). The parameter $\theta$, defined in Section 1, satisfies $\theta = \theta(P^T)$. The estimate of $\theta$ is easy in some cases:

(a) if there exists $\nu : X \mapsto [0, +\infty)$, $\nu \not\equiv 0$, such that $\nu(x) \leq C\nu(y)$ for any $x, y \in X$ (and some constant $C > 0$) and $\nu(x) \geq \sum_y \nu(y) p(y,x)$ [resp. $\nu(x) \leq \sum_y \nu(y) p(y,x)$], then $\theta \geq 1$ (resp. $\theta \leq 1$);
(b) if $\lim_{n \to +\infty} (\#\{y : |y| \leq n\})^{1/n} = 1$, then $\theta \leq 1$.

Note that in general $\rho \leq 1$ and $\rho \leq \theta$. Using the above result (a), given a graph $X$, if $P$ is *strongly reversible* (see [6]), then the amenability of $X$ implies $\rho = \theta = 1$, while nonamenability implies $\rho < \theta = 1$. Roughly speaking, amenable graphs are graphs where the boundary of a finite set may be arbitrarily small if compared to the size of the set itself (see [16], page 112). Hence, $\rho = \theta = 1$ in the case of the simple random walk on $\mathbb{Z}^d$, while there are examples where $\theta \neq 1$ and $\theta > \rho$, for instance, on homogeneous trees $\mathbb{T}_{n+1}$ [see Example 5.1 with $p \neq 1/(n+1)$].

2.2. *Configuration space.* Following [10], fix a reference vertex $x_0 \in X$ and denote by $|x|$ the graph distance between $x$ and $x_0$. Define a strictly positive function $\alpha \colon X \to \mathbb{R}^+$ by $\alpha(x) = M^{-|x|}$, where $M > (D-1)^2$. By this choice of $M$, for any $z \geq 1/2$, $\sum_x \alpha(x)^z < +\infty$ and

$$\sum_y q(x,y)\alpha(y) \leq M\alpha(x) \tag{2.1}$$

for any substochastic matrix $Q$. Given $\eta \colon X \to \mathbb{N}$, define $\|\eta\| := \sum_x \eta(x)\alpha(x)$. The configuration space is $\Omega := \{\eta \in \mathbb{N}^X \text{ such that } \|\eta\| < +\infty\}$, while $\Omega_\Lambda := \{\eta \colon \mathbb{N}^\Lambda \text{ such that } \|\eta\| < +\infty\}$. Note that the *finite* configurations, that is, the configurations $\eta \in \Omega$ such that $\sum_x \eta(x) < +\infty$, are dense in $\Omega$ with this norm; moreover, the Borel $\sigma$-algebra induced by the norm is the same as the one induced by the product topology. We introduce the usual partial order on $\Omega$, that is, $\xi \leq \eta$ if $\xi(x) \leq \eta(x)$ for any $x \in X$. We denote by $\mathbf{0}$ the configuration identically equal to 0, by $\mathbf{1}$ the configuration identically equal to 1 and by $\delta_x$ the configuration which is equal to 0 at any site but $x$, where it equals 1. We say that a function $f \colon \Omega \to \mathbb{R}$ is *nondecreasing* if $\xi \leq \eta$ implies $f(\xi) \leq f(\eta)$. Given $\mu, \nu$ probability measures on $\Omega$, we say that $\nu$ *stochastically dominates* $\mu$ and we write $\mu \leq \nu$ if for any nondecreasing function $f$ we have $\mu(f) \leq \nu(f)$ [where $\mu(f) = \int_\Omega f \, d\mu$].

From (2.1) we derive a useful bound on the transition kernel of the continuous time random walk associated with $P$ and with jump rate $\lambda > 0$. Indeed, let

$$p_t^\lambda(x,y) := e^{-\lambda t} \sum_{n=0}^{+\infty} \frac{(\lambda t)^n}{n!} p^{(n)}(x,y),$$

then the iteration of (2.1) gives

$$\sum_y p_t^\lambda(x,y)\alpha(y) \leq e^{\lambda t(M-1)}\alpha(x). \tag{2.2}$$

For $\Lambda \subset X$, we will also denote by

$$p_{t,\Lambda}^\lambda(x,y) := e^{-\lambda t} \sum_{n=0}^{+\infty} \frac{(\lambda t)^n}{n!} p_\Lambda^{(n)}(x,y). \tag{2.3}$$

Clearly, $p_{t,\Lambda}^\lambda(x,y) \leq p_t^\lambda(x,y)$, hence, the bound in (2.2) holds for these "restricted" kernels as well.

2.3. *Dynamics.* Denote by $\text{Lip}(\Omega)$ the set of the Lipschitz functions on $\Omega$, and given $f \in \text{Lip}(\Omega)$, let $L(f)$ be its Lipschitz constant. For any $f \colon \Omega \to \mathbb{R}$ and $x \in X$, define $(\partial_x^- f)(\eta) = \mathbb{1}_{[1,+\infty)}(\eta(x))[f(\eta - \delta_x) - f(\eta)]$ and $(\partial_x^+ f)(\eta) = f(\eta + \delta_x) - f(\eta)$. Note that $|(\partial_x^\pm f)(\eta)| \leq L(f)\alpha(x)$ for any $f \in \text{Lip}(\Omega)$, $x \in X$.

Fix a nonincreasing function $c: \mathbb{N} \to \mathbb{R}^+$, a transition matrix $P$ on $X$ and define $\mathcal{L}: \text{Lip}(\Omega) \to \mathbb{R}^\Omega$ by

$$(2.4) \quad (\mathcal{L}f)(\eta) := \sum_x \eta(x) \left[ (\partial_x^- f)(\eta) + \sum_y c(\eta(y)) p(x,y)(\partial_y^+ f)(\eta) \right].$$

It easy to check that $|(\mathcal{L}f)(\eta)| \leq L(f)[(c(0)M+1)\|\eta\|]$, hence, $\mathcal{L}$ is a well-defined operator on $\text{Lip}(\Omega)$.

## 3. Construction of the process.

The main result of this section concerns the construction of a process having generator $\mathcal{L}$ given by (2.4). It is a standard fact of the theory of countable state continuous time Markov chains that there exists a unique Markov process $\{\eta_t\}_{t \geq 0}$ with generator given by (2.4) starting from any finite $\eta \in \Omega$. The extension of this construction to more general configurations requires more sophisticated techniques. Our efforts in this direction may be summarized in the following proposition whose proof is the consequence of several intermediate steps.

PROPOSITION 3.1. *There exists a unique semigroup $\{S_t\}_{t \geq 0}$ of operators, $S_t: \text{Lip}(\Omega) \to \text{Lip}(\Omega)$ such that:*

(i) $(S_t f)(\eta) = \mathbb{E}^\eta[f(\eta_t)]$ *for* $f \in \text{Lip}(\Omega)$ *and $\eta$ finite.*
(ii) $L(S_t f) \leq L(f) \exp(c(0)Mt)$ *for* $f \in \text{Lip}(\Omega)$.
(iii) $(S_t f)(\eta) = f(\eta) + \int_0^t (\mathcal{L} S_u f)(\eta)\, du$ *for* $f \in \text{Lip}(\Omega)$ *and* $\eta \in \Omega$.
(iv) *Let $\mu$ be a probability measure on $\Omega$ such that $\mu[\|\eta\|]$ is finite. Then $\mu$ is invariant for $\{S_t\}_{t \geq 0}$ [i.e., for any $t \geq 0$, $\mu[S_t f] = \mu[f]$ for any $f \in \text{Lip}(\Omega)$] if and only if $\mu[\mathcal{L}f] = 0$ for any $f \in \text{Lip}(\Omega)$.*

Given Proposition 3.1, in order to define the process $\{\eta_t\}_{t \geq 0}$ starting from any $\eta \in \Omega$, according to [10] (see discussion after Theorem 1.4), one shows that $\{S_t\}_{t \geq 0}$ can be extended to any measurable function $f$ on $\Omega$ which satisfies either $f \geq 0$ or $|f(\eta)| \leq C(1 + \|\eta\|)$ for some constant $C$. Thus, it identifies a unique Markov process on $\Omega$ (the RBRW) which we still denote by $\{\eta_t\}_{t \geq 0}$.

We start by constructing the process on a finite subset $\Lambda \subset X$. We need an auxiliary process defined on $\Omega$. Fix $\gamma \geq 0$, $c: \mathbb{N} \to \mathbb{R}^+$, $k \in \mathbb{N}$ and define $\mathcal{G}_\Lambda: \text{Lip}(\Omega) \to \mathbb{R}^\Omega$ by

$$(3.1) \quad (\mathcal{G}_\Lambda f)(\eta) := \sum_x \left[ \gamma \mathbb{1}_\Lambda(x)(\eta(x) - k)^+ (\partial_x^- f)(\eta) + \eta(x) \sum_y p_\Lambda(x,y) c(\eta(y))(\partial_y^+ f)(\eta) \right].$$

It is obvious that $\mathcal{G}_\Lambda$ generates a Markov process $\{\eta_t\}_{t\geq 0}$ defined on $\Omega_\Lambda$. Clearly, this process may be thought as a process on $\Omega$ where the particles outside $\Lambda$ are "frozen" in the initial state. Furthermore, if this process starts from $\eta_0 \in \Omega$ such that $\eta_0 \geq k\mathbf{1}$, then obviously $\eta_t \geq k\mathbf{1}$ for any $t \geq 0$; in this case we say that there are $k$ *immortal particles* per site.

3.1. *BRW with immortal particles.* The estimate of the first and the second moments of the process generated by (3.1) will follow from a coupling with the BRW with $k$ immortal particles, that is, the process where $c(\cdot) \equiv \lambda$. Hence, we take a technical detour and study this particular process [or, equivalently, the BRW with constant immigration rate; see (3.4)].

Although in this section we are treating only finite sets, the following lemma is needed also in the countable case (see Remark 3.11).

LEMMA 3.2. *Let $Q$ be a (possibly infinite) matrix with $q(x,y) \in [0,1]$ and finite norm and define $q_t(x,y) := e^{-\lambda t}\sum_{n=0}^{\infty}(\lambda t)^n q^{(n)}(x,y)/n!$. Let $f:[0,+\infty) \to l^\infty(X)$ be such that $\lim_{t\to\infty} f(t) = v$. Given the system of linear differential equations,*

$$(3.2) \quad \begin{cases} \dot{u}(t,x) = \lambda\bigg(\sum_y q(x,y)u(t,y) - u(t,x)\bigg) + \beta u(t,x) + f(t,x) & \forall x, \\ u(0,\cdot) = \varphi(\cdot), \end{cases}$$

*where $\varphi \in l^\infty(X)$ and $\beta > \lambda(1 + \theta(Q))$, the corresponding solution satisfies $\lim_{t\to\infty} u(t,x) = ((\lambda - \beta)\mathbb{I} - \lambda Q)^{-1}v$ for any $x \in X$ and $\varphi \in l^\infty(X)$.*

PROOF. The proof is standard and we just sketch it. One can solve the system by considering the (stronger) Cauchy problem in $l^\infty(X)$

$$\begin{cases} \dot{u}(t) = -Au(t) + f(t), \\ u(0) = \varphi, \end{cases}$$

where $A = (\lambda - \beta)\mathbb{I} - \lambda Q$. By our hypotheses we have that $\mathbf{Re}(\sigma(A)) \geq \epsilon > 0$, hence, $\|e^{-At}\| \xrightarrow{t\to\infty} 0$, $\|\int_0^t e^{-As}\,ds\| \leq \int_0^{+\infty}\|e^{-As}\|\,ds < +\infty$ for any $t \in [0,+\infty]$ and $A$ is invertible. The solution is given by the well-known formula $u(t) = e^{-At}\varphi + \int_0^t e^{-A(t-s)}f(s)\,ds$. The first term tends to zero, while the second one can be written as $\int_0^{t_0} e^{-A(t-s)}f(s)\,ds + \int_{t_0}^t e^{-A(t-s)}f(s)\,ds$ and the claim follows choosing $t_0$ such that for $t \geq t_0$, $\|f(t) - v\|_\infty$ is sufficiently small. □

LEMMA 3.3. *Let $\Lambda \subset X$ be finite. Fix $\gamma \geq 0$ and $k \in \mathbb{N}$. Consider $\mathcal{G}_\Lambda$ defined in (3.1) with $c(\cdot) \equiv \lambda > 0$. Let $\{\eta_t\}_{t\geq 0}$ be the process generated by $\mathcal{G}_\Lambda$ starting from $\eta$ bounded, $\eta \geq k\mathbf{1}$. Moreover, if $\lambda < \gamma/\theta$, then there exists* *two nonnegative constants $U_{1,\Lambda}(k,\lambda,\gamma)$ and $U_{2,\Lambda}(k,\lambda,\gamma)$ such that, for any $x \in \Lambda$, we have that*

$$(3.3) \quad \lim_{t \to \infty} \mathbb{E}^\eta[\eta_t(x)] \leq U_{1,\Lambda}(k,\lambda,\gamma), \qquad \lim_{t \to \infty} \mathbb{E}^\eta[(\eta_t(x))^2] \leq U_{2,\Lambda}(k,\lambda,\gamma),$$

*where the limits are attained uniformly with respect to $x$.*

PROOF. Define $\{\xi_t\}_{t \geq 0}$ as $\xi_t := \eta_t - k\mathbf{1}$. This process is a Markov process (viz., it is the branching random walk with constant immigration rate $\lambda k$) and its generator is

$$(3.4) \quad (\mathcal{H}_\Lambda g)(\xi) := \sum_{x \in \Lambda} \left[ \gamma \xi(x)(\partial_x^- g)(\xi) + \lambda(\xi(x) + k) \sum_y p_\Lambda(x,y)(\partial_y^+ g)(\xi) \right].$$

Obviously for any $\eta \in \Omega$ such that $\eta \geq k\mathbf{1}$, we have

$$(3.5) \quad \begin{aligned} \mathbb{E}^\eta[\eta_t(x)] &= \mathbb{E}^{\eta - k\mathbf{1}}[\xi_t(x)] + k, \\ \mathbb{E}^\eta[(\eta_t(x))^2] &= \mathbb{E}^{\eta - k\mathbf{1}}[(\xi_t(x))^2] + 2k\mathbb{E}^{\eta - k\mathbf{1}}[\xi_t(x)] + k^2. \end{aligned}$$

Choose $\xi \in \Omega$ and let $m(t,x) := \mathbb{E}^\xi[\xi_t(x)]$, for any $x \in \Lambda$; by basic semigroup properties we have that $\frac{d}{dt} m(t,x) = \mathbb{E}^\xi[(\mathcal{H}_\Lambda \pi_x)(\xi_t)]$ (where $\pi_x$ is the projection on the $x$ coordinate). By computing explicitly $\mathcal{H}_\Lambda \pi_x$, we obtain that $m$ satisfies the system (3.2) with $Q = P^T$, $f(t,x) = k\lambda \sum_y p_\Lambda(y,x)$ and $\beta = \lambda - \gamma$. The claim follows from Lemma 3.2.

To prove the second moment assertion, consider $C(t,x,y) := \mathbb{E}^\xi[\xi_t(x)\xi_t(y)]$, for any $x,y \in \Lambda$. Using the same arguments as before, we obtain that $C$ is the solution following system of linear differential equations:

$$(3.6) \quad \begin{cases} \dfrac{d}{dt} C(t,x,y) = -AC(t,x,y) + f(t,x,y) & \forall x,y \in \Lambda, \\ C(0,x,y) = \xi(x)\xi(y), \end{cases}$$

where

$$A = 2(\gamma \mathbb{I} - \lambda B),$$

$$B((x,y),(x_1,y_1)) = \delta_y(y_1)\frac{p_\Lambda(x_1,x)}{2} + \delta_x(x_1)\frac{p_\Lambda(y_1,y)}{2},$$

$$f(t,x,y) = \lambda k \left( m(t,x) \sum_z p_\Lambda(z,y) + m(t,y) \sum_z p_\Lambda(z,x) \right)$$
$$+ \delta_x(y)\lambda \left( k \sum_z p_\Lambda(z,x) + \sum_z p_\Lambda(z,x)m(t,z) \right)$$
$$+ \delta_x(y)\gamma m(t,x).$$

The system (3.6) is formally equivalent to the one in (3.2) with $X \times X$ in the place of $X$. The results just obtained for $m$ ensure that $f$ satisfies the assumptions in Lemma 3.2. Moreover, $\theta(B) \le \theta$, since $b^{(n)}((x,y),(x_1,y_1)) = \frac{1}{2^n}\sum_{k=0}^{n}\binom{n}{k}p^{(k)}(x_1,x)p^{(n-k)}(y_1,y)$ and the claim follows. $\square$

REMARK 3.4. It is known that, given the equations (3.2) with $x \in \Lambda$, where $\Lambda$ is finite, an explicit expression of the solution, for any $\varphi \in \Omega$, is

$$u(t,x) = e^{\beta t}\sum_y q_t(x,y)\varphi(y) + \sum_y \int_0^t e^{\beta(t-s)}q_{t-s}(x,y)f(s,y)\,ds.$$

This formula represents the solution also for infinite $\Lambda$ under mild assumptions on $f$: suppose, for instance, that $f \ge 0$, $\sum_x f(t,x)\alpha(x) < +\infty$ for some $t \ge 0$ and $\partial_t f(t,x)$ is bounded on any compact set, uniformly with respect to $x \in \Lambda$. Moreover, if we consider two families $\{Q_\Lambda\}_\Lambda$ and $\{f_\Lambda\}_\Lambda$ which are nondecreasing with respect to $\Lambda$ and we denote by $\{u_\Lambda\}_\Lambda$ the corresponding solutions, then we have that $u_\Lambda \uparrow u_X$ as $\Lambda \uparrow X$ [hence, one may replace the upper bounds in (3.3) with $U_{1,X}$ and $U_{2,X}$ which clearly are uniform with respect to $\Lambda$].

We note that the bound (3.3) can be obtained starting from any $\eta \in \Omega$ (not necessarily bounded) by using similar computations, if $\lambda < \gamma/\|P^T\| \le \gamma/\theta$. It is easy to show that $\|P^T\| = 1$, for instance, for any symmetric random walk.

3.2. *The finite volume process.* The following proposition shows how to construct a monotone coupling of different processes generated by (3.1).

PROPOSITION 3.5. *Fix $N \in \mathbb{N}$ and let $\Lambda_1 \subset \cdots \subset \Lambda_N \subset X$ be finite subsets. Fix $k_1,\ldots,k_N \in \mathbb{N}$, $\gamma_1,\ldots,\gamma_N \in [0,+\infty)$ and let $c_1,\ldots,c_N : \mathbb{N} \to [0,+\infty)$ be nonincreasing functions. Then for any fixed $(\eta_{0,1},\ldots,\eta_{0,N}) \in \Omega^N$ such that $\eta_{0,h} \ge k_h\mathbf{1}$ for any $h \in \{1,\ldots,N\}$, there exists a Markov process $\{(\eta_{t,1},\ldots,\eta_{t,N})\}_{t\ge 0}$ on $\Omega^N$ such that, for any $h \in \{1,\ldots,N\}$, the semigroup associated with the process $\{\eta_{t,h}\}_{t\ge 0}$ has generator $\mathcal{G}_{\Lambda_h}$. Furthermore, assume that $k_h \le k_{h+1}$, $\gamma_h \ge \gamma_{h+1}$, $c_h(k_{h+1}+n) \le c_{h+1}(k_{h+1}+n)$ for any $n \in \mathbb{N}$ and $\eta_{0,h} \le \eta_{0,h+1}$ for any $h \in \{1,\ldots,N-1\}$. Then $\eta_{t,1} \le \cdots \le \eta_{t,N}$ for any $t \ge 0$.*

PROOF. It is enough to consider the processes on $\Omega_{\Lambda_N}$. Choose $(\eta_{0,1},\ldots,\eta_{0,N})$ in $(\Omega_{\Lambda_N})^N$ such that $\eta_{0,h} \ge k_h\mathbf{1}$ for any $h \in \{1,\ldots,N\}$ as the initial configurations. For any $x \in \Lambda_N$, let $A(x) := \max\{\eta_{0,1}(x),\ldots,\eta_{0,N}(x)\}$, $\bar\gamma := \max\{\gamma_1,\ldots,\gamma_N\}$, $\bar c := \max\{c_1(0),\ldots,c_N(0)\}$. Choose an independent family of exponential clocks, two per site $x \in \Lambda_N$: one of parameter $\bar\gamma A(x)$ which controls the deaths and one of parameter $\bar c A(x)$ which controls births. Define $(\eta_{t,1},\ldots,\eta_{t,N}) := (\eta_{0,1},\ldots,\eta_{0,N})$ for any $t < \tau$, where $\tau$ is the time of the first ring of the collection of clocks. Assume that the clock which rings first is at site $x$:

- If the clock is a death clock, then for any $z \neq x$ put $\eta_{\tau,h}(z) := \eta_{\tau-,h}(z)$, pick a uniform $U$ in the interval $(0,1)$ and define $\eta_{\tau,h}(x) := \eta_{\tau-,h}(x) - 1$ for any $h$ such that $U \leq (\gamma_h(\eta_{\tau-,h}(x) - k_h)^+)/\bar{\gamma} A(x)$ and $\eta_{\tau,h}(x) := \eta_{\tau-,h}(x)$ otherwise. Finally, restart the procedure from $(\eta_{\tau,1}, \ldots, \eta_{\tau,N})$.
- If the clock is a birth clock, then for any $z \not\sim x$ put $\eta_{\tau,h}(z) := \eta_{\tau-,h}(z)$. Choose at random, accordingly to the transition matrix $P$, a site $y$ among the neighbors of $x$. Now pick a uniform $V$ in $(0,1)$ and define $\eta_{\tau,h}(y) := \eta_{\tau-,h}(y) + 1$ for any given $h$ such that $V$ is not larger than $\eta_{\tau-,h}(x) \times p_{\Lambda_h}(x,y) c_h(\eta_{\tau-,h}(y))/(\bar{c} A(x) p(x,y))$ and $\eta_{\tau,h}(y) = \eta_{\tau-,h}(y)$ otherwise. Finally, restart the procedure from $(\eta_{\tau,1}, \ldots, \eta_{\tau,N})$.

It is a simple exercise to check that this construction leads to the desired coupling. □

In the remaining part of this section we prove some basic bounds on the semigroup $\{S_{t,\Lambda}\}_{t \geq 0}$ generated by $\mathcal{G}_\Lambda$. We need these bounds to extend the construction of the process to an infinite $\Lambda \subset X$. The next result shows that the semigroup $\{S_{t,\Lambda}\}_{t \geq 0}$ generated by $\mathcal{G}_\Lambda$ maps $\mathrm{Lip}(\Omega)$ into itself. The proof follows closely the proof of Lemma 2.1 in [10].

LEMMA 3.6. *Let $\Lambda \subset X$ be finite and $\{S_{t,\Lambda}\}_{t \geq 0}$ be the semigroup generated by $\mathcal{G}_\Lambda$. Then for any $f \in \mathrm{Lip}(\Omega)$,*

$$L(S_{t,\Lambda} f) \leq L(f) \exp(c(0) M t).$$

PROOF. Take $\xi, \zeta \in \Omega$ $(\xi, \zeta \geq k\mathbf{1})$ and consider the monotone coupling $\{(\eta_t^1, \eta_t^2, \eta_t^3, \eta_t^4)\}_{t \geq 0}$ of Proposition 3.5 such that $(\eta_0^1, \eta_0^2, \eta_0^3, \eta_0^4) := (\xi \wedge \zeta, \xi, \zeta, \xi \vee \zeta)$. This means that $\eta_t^1 \leq \eta_t^2, \eta_t^3 \leq \eta_t^4$ for any $t \geq 0$. Therefore,

$$|(S_{t,\Lambda} f)(\xi) - (S_{t,\Lambda} f)(\zeta)| = |\mathbb{E}[f(\eta_t^2) - f(\eta_t^3)]|$$
$$\leq \mathbb{E}[|f(\eta_t^2) - f(\eta_t^3)|] \leq L(f) \mathbb{E}[\|\eta_t^2 - \eta_t^3\|].$$

To bound this last term notice that by monotonicity

$$\|\eta_t^2 - \eta_t^3\| \leq \sum_x \alpha(x)(\eta_t^4(x) - \eta_t^1(x)).$$

Furthermore, for any $x \in \Lambda$, we claim that

$$(3.7) \qquad \frac{d}{dt} \mathbb{E}[\eta_t^4(x) - \eta_t^1(x)] \leq c(0) \sum_y p_\Lambda(y,x) \mathbb{E}[\eta_t^4(y) - \eta_t^1(y)],$$

which implies, by (2.1), that

$$\frac{d}{dt} \mathbb{E}[\|\eta_t^4 - \eta_t^1\|] \leq c(0) M \mathbb{E}[\|\eta_t^4 - \eta_t^1\|].$$

This gives
$$\mathbb{E}[\|\eta_t^4 - \eta_t^1\|] \leq \mathbb{E}[\|\eta_0^4 - \eta_0^1\|]\exp(c(0)Mt)$$
and the proof is complete. To obtain (3.7), use the generator of the coupled process or better notice (see the proof of Proposition 3.5) that the rate of the transition $\eta_t^4(x) - \eta_t^1(x) \to \eta_t^4(x) - \eta_t^1(x) + 1$ is
$$c(\eta_{t-}^4(x))\sum_y p_\Lambda(y,x)(\eta_{t-}^4(y) - \eta_{t-}^1(y)). \qquad \square$$

The following result is a simple consequence of Lemma 3.6. The proof is the same as the one of Corollary 2.5 in [10], hence, we omit it.

COROLLARY 3.7. *Let $\Lambda \subset X$ be finite and $\{S_{t,\Lambda}\}_{t \geq 0}$ be the semigroup generated by $\mathcal{G}_\Lambda$. For any $f : \Omega \to \mathbb{R}$ such that $|f(\eta)| \leq C_f \|\eta\|$ for all $\eta \geq k\mathbf{1}$ and for some constant $C_f > 0$, we have that*
$$|(S_{t,\Lambda}f)(\eta)| \leq C_f\|\eta\|\exp(c(0)Mt).$$

The next two results are the analogs of Lemma 2.6 and Lemma 2.7 in [10].

LEMMA 3.8. *Let $\Lambda \subset \Lambda' \subset X$ be finite subsets. Fix $\gamma \geq 0$ and $k \in \mathbb{N}$. Let $\mathcal{G}_\Lambda, \mathcal{G}_{\Lambda'}$ be defined by (3.1). Then for any $f \in \mathrm{Lip}(\Omega)$, $\eta \geq k\mathbf{1}$,*
$$|(\mathcal{G}_{\Lambda'}f)(\eta) - (\mathcal{G}_\Lambda f)(\eta)| \leq L(f)(\gamma + c(0)M)\sum_x \mathbb{1}_{\Lambda'\setminus\Lambda^0}(x)\eta(x)\alpha(x).$$

PROOF. The proof can be obtained by direct computation. $\square$

LEMMA 3.9. *Let $\Lambda \subset \Lambda' \subset X$ be two finite subsets. Fix $\gamma \geq 0$ and $k \in \mathbb{N}$. Consider the semigroups $\{S_{t,\Lambda}\}_{t\geq 0}$ and $\{S_{t,\Lambda'}\}_{t\geq 0}$ associated with the generators $\mathcal{G}_\Lambda, \mathcal{G}_{\Lambda'}$ defined by (3.1). Then for any $f \in \mathrm{Lip}(\Omega)$ and $\eta \geq k\mathbf{1}$,*
$$|(S_{t,\Lambda'}f)(\eta) - (S_{t,\Lambda}f)(\eta)|$$
$$\leq L(f)(\gamma + c(0)M)e^{c(0)Mt}\sum_{x,y}\alpha(x)\mathbb{1}_{\Lambda'\setminus\Lambda^0}(x)\eta(y)\int_0^t p_{u,\Lambda'}^{c(0)}(y,x)\,du.$$

PROOF. Note that
$$(S_{t,\Lambda'}f)(\eta) - (S_{t,\Lambda}f)(\eta) = \int_0^t (S_{u,\Lambda'}(\mathcal{G}_{\Lambda'} - \mathcal{G}_\Lambda)S_{t-u,\Lambda}f)(\eta)\,du.$$
By Lemma 3.8,
$$|((\mathcal{G}_{\Lambda'} - \mathcal{G}_\Lambda)S_{t-u,\Lambda}f)(\eta)| \leq L(S_{t-u,\Lambda}f)(\gamma + c(0)M)\sum_x \mathbb{1}_{\Lambda'\setminus\Lambda^0}(x)\eta(x)\alpha(x).$$

By Lemma 3.6,
$$L(S_{t-u,\Lambda}f) \leq L(f)e^{c(0)M(t-u)}.$$

Using this last estimate and the positivity of $S_{u,\Lambda'}$, we get

(3.8)
$$|(S_{u,\Lambda'}(\mathcal{G}_{\Lambda'} - \mathcal{G}_{\Lambda})S_{t-u,\Lambda}f)(\eta)|$$
$$\leq L(f)e^{c(0)M(t-u)}(\gamma + c(0)M)\sum_x \alpha(x)\mathbb{1}_{\Lambda'\setminus\Lambda^0}(x)S_{u,\Lambda'}(\pi_x)(\eta).$$

By Proposition 3.5, $(S_{u,\Lambda'}\pi_x)(\eta) \leq \mathbb{E}^\eta_{\Lambda'}[\eta_u(x)]$, where $\{\eta_t\}_{t\geq 0}$ is the process generated by (3.1) with $\gamma = 0$, $c(\cdot) \equiv c(0)$ and $k = 0$. By Remark 3.4, we know that
$$\mathbb{E}^\eta_{\Lambda'}[\eta_u(x)] = e^{c(0)u}\sum_y p^{c(0)}_{u,\Lambda'}(y,x)\eta(y).$$

Plugging this bound in (3.8), we get
$$|(S_{u,\Lambda'}(\mathcal{G}_{\Lambda'} - \mathcal{G}_{\Lambda})S_{t-u,\Lambda}f)(\eta)|$$
$$\leq L(f)e^{c(0)Mt}(\gamma + c(0)M)\sum_{x,y} \alpha(x)\mathbb{1}_{\Lambda'\setminus\Lambda^0}(x)p^{c(0)}_{u,\Lambda'}(y,x)\eta(y),$$

which concludes the proof. □

3.3. *Finite volume approximation.* Following [10], we construct the process on $X$ as a limit of processes defined on $\Lambda$ finite. For any $n \in \mathbb{N}$, define $\Lambda_n := B(x_0, n)$, that is, the ball of radius $n$ and center $x_0$.

PROPOSITION 3.10. *Fix $\gamma \geq 0$ and $k \in \mathbb{N}$. For any $n \in \mathbb{N}$, consider the semigroups $\{S_{t,\Lambda_n}\}_{t\geq 0}$ generated by $\mathcal{G}_{\Lambda_n}$ defined in (3.1). For any fixed $t \geq 0$, $f \in \operatorname{Lip}(\Omega)$ and $\eta \in \Omega$, $\eta \geq k\mathbf{1}$, the sequence $\{S_{t,\Lambda_n} : n \in \mathbb{N}\}$ is a Cauchy sequence.*

PROOF. Assume that $m \leq n$, then by Lemma 3.9,
$$|(S_{t,\Lambda_n}f)(\eta) - (S_{t,\Lambda_m}f)(\eta)|$$
$$\leq L(f)(\gamma + c(0)M)e^{c(0)Mt}\sum_{x,y}\alpha(x)\mathbb{1}_{\Lambda_n\setminus\Lambda^0_m}(x)\eta(y)\int_0^t p^{c(0)}_{u,\Lambda_n}(y,x)\,du.$$

We have to show that
$$\sum_y p^{c(0)}_{u,\Lambda_n}(y,x)\eta(y)$$

can be dominated uniformly in $n \in \mathbb{N}$ by a function $\phi(x,u) \in L^1(X \times [0,t])$, where the measure on $X$ is $\alpha(\cdot)$. The result follows by dominated convergence since $\lim_{m,n\to+\infty} \mathbb{1}_{\Lambda_n\setminus\Lambda_m^0}(x) \le \lim_{m\to+\infty} \mathbb{1}_{(\Lambda_m^0)^\complement}(x) = 0$. We claim that we can take $\phi(x,u) = \sum_y \eta(y) p_u^{c(0)}(y,x)$, indeed, $p_{u,\Lambda_n}^{c(0)}(y,x) \le p_u^{c(0)}(y,x)$ and by (2.2), we have

$$\sum_x \alpha(x) \sum_y p_u^{c(0)}(y,x)\eta(y) \le e^{c(0)(M-1)u} \sum_y \eta(y)\alpha(y)$$
$$= e^{c(0)(M-1)u}\|\eta\| \in L^1([0,t]). \qquad \square$$

The proposition above allows us to define for any $t \ge 0$, $f \in \text{Lip}(\Omega)$ and $\eta \in \Omega$, $\eta \ge k\mathbf{1}$:

$$(S_t f)(\eta) := \lim_{n\to+\infty} (S_{t,\Lambda_n} f)(\eta).$$

REMARK 3.11. It easy to show that with this definition we can drop the hypothesis that $\Lambda$ is finite (take $\Lambda \uparrow X$) in Proposition 3.5, Lemma 3.6, Corollary 3.7, Lemma 3.8 and Lemma 3.9. The same can be done in Lemma 3.3, since one proves that $\mathbb{E}^\eta[\eta_t(x)] = \lim_{\Lambda\uparrow X} m(t,x) + k$ and $\mathbb{E}^\eta[\eta_t(x)\eta_t(y)] = \lim_{\Lambda\uparrow X} C(t,x,y)$, where the limit functions satisfy the corresponding differential systems. Note that, in particular, the latter is not obvious, because $\eta \mapsto \eta(x)\eta(y) \notin \text{Lip}(\Omega)$. Moreover, the process generated by (2.4) is monotone as a consequence of Proposition 3.5.

PROPOSITION 3.12. *For any $t \ge 0$, $f \in \text{Lip}(\Omega)$ and $\eta \in \Omega$, $\eta \ge k\mathbf{1}$, define*

$$(S_t f)(\eta) := \lim_{n\to+\infty} (S_{t,\Lambda_n} f)(\eta):$$

1. $\{S_t\}_{t \ge 0}$ *is a semigroup.*
2. *For all $f \in \text{Lip}(\Omega)$, $\eta \in \Omega$,*

$$(S_t f)(\eta) = f(\eta) + \int_0^t (\mathcal{G} S_u f)(\eta)\,du.$$

PROOF. These properties can be proven exactly as in [10] (page 451 and Lemma 2.12) by using Lemma 3.6, Corollary 3.7, Lemma 3.3, Lemma 3.8 and Lemma 3.9 instead of Lemma 2.1, Corollary 2.5, Lemma 2.6 and Lemma 2.7, respectively. $\square$

Among the properties of the semigroup $\{S_t\}_{t\ge 0}$ which can be proven, we state the one which we need in the next section.

PROPOSITION 3.13. *Let $\mu$ be a probability measure on $\Omega$ such that $\mu[\|\eta\|]$ is finite. Then $\mu$ is invariant for $\{S_t\}_{t\ge 0}$ (i.e., for any $t \ge 0$, $\mu[S_t f] = \mu[f]$ for any $f \in \text{Lip}(\Omega)$) if and only if $\mu[\mathcal{G}f] = 0$ for any $f \in \text{Lip}(\Omega)$.*

PROOF. See the proof of Corollary 2.17 in [10]. □

PROOF OF PROPOSITION 3.1. It is easy to show that (i) holds. The claim (ii) follows from Lemma 3.6 and Remark 3.11, while Propositions 3.12 and 3.13 imply (iii) and (iv), respectively. □

**4. Ecological equilibrium and invariant measure.** In this section we study the behavior of the RBRW constructed in Section 3. In particular, Proposition 1.1 and Theorem 1.2 are proven (see Proposition 4.1 and Theorem 4.3 below). The main tool is the coupling between this monotone process and suitable contact and BRW processes.

PROPOSITION 4.1. *Let $\{\eta_t\}_{t\geq 0}$ the RBRW generated by (2.4):*

(i) *If $c(0) \leq 1/\rho$ [resp. $c(0) < 1/\rho$], then $\lim_{t\to+\infty} \eta_t(x) = 0$ a.s. (resp. $\lim_{t\to+\infty} \mathbb{E}^\eta[\eta_t(x)] = 0$) for any finite $\eta \in \Omega$, $x \in X$;*

(ii) *If $c(0) > \lambda_{\mathrm{CP}}$, then $\lim_{t\to+\infty} \mathbb{E}^\eta[\eta_t(x)] > 0$ for any $\eta \in \Omega \setminus \{\mathbf{0}\}$, $x \in X$ and $\mathbb{P}^\eta(\limsup_{t\to\infty} \eta_t(x) > 0) > 0$;*

(iii) *If $c(+\infty) > 1/\rho$, then $\lim_{t\to+\infty} \mathbb{E}^\eta[\eta_t(x)] = +\infty$ for any $\eta \in \Omega \setminus \{\mathbf{0}\}$;*

(iv) *If $c(+\infty) < 1/\theta$, then $\limsup_{t\to+\infty} \mathbb{E}^\eta[\eta_t(x)] < +\infty$ uniformly for any bounded $\eta \in \Omega$, $x \in X$.*

PROOF. Recall that Proposition 3.5 and Lemma 3.3 hold for $\Lambda = X$ (see Remark 3.11):

(i) By Proposition 3.5, we can couple the process with a branching random walk $\{\zeta_t\}_{t\geq 0}$ starting from $\eta$ with birth rate $c(0)$ such that $\eta_t \leq \zeta_t$. The first part of the claim follows by noting that $\zeta_t$ dies out almost surely (see [1], Theorem 3.1). As for the second part, the assertion is a consequence of Lemma 3.3 and Remark 3.11 [since $k = 0$, one can choose $U_{1,X}(0, \lambda, 1) = 0$].

(ii) By Proposition 3.5, there exists a supercritical site-breeding CP $\{\zeta_t\}_{t\geq 0}$ starting from $\eta \wedge \mathbf{1}$ with birth rate $c(0) > \lambda_{\mathrm{CP}}$ and such that $\zeta_t \leq \eta_t$. Theorem 4.8 of Chapter VI in [8] yields to the conclusion.

(iii) By Proposition 3.5, we can couple $\{\eta_t\}_{t\geq 0}$ with a branching random walk $\{\zeta_t\}_{t\geq 0}$ starting from $\eta$ with birth rate $c(0) > 1$ such that $\zeta_t \leq \eta_t$ and the claim follows.

(iv) In this case there exists $\bar{k} \in \mathbb{N}$ such that $c(\bar{k}) < 1/\theta$. By Proposition 3.5 and Remark 3.11, there exists a process $\{\zeta_t\}_{t\geq 0}$ generated by (3.1) with $k = \bar{k}$, $\gamma = 1$, birth rate $c(\bar{k})$, such that $\zeta_t \geq \eta_t$. By Lemma 3.3 and Remark 3.11, we have

$$\limsup_{t\to\infty} \mathbb{E}^{\eta\vee \bar{k}\mathbf{1}}[\eta_t(x)] \leq \lim_{t\to\infty} \mathbb{E}^{\eta\vee \bar{k}\mathbf{1}}[\zeta_t(x)] \leq U_{1,X}(\bar{k}, c(\bar{k}), 1).$$

□

REMARK 4.2. The condition $\limsup_{t \to +\infty} \mathbb{E}^\eta[\eta_t(x)] < +\infty$ implies that $\mathbb{P}^\eta(\lim_{t \to +\infty} \eta_t(x) = +\infty) = 0$, but $\eta_t(x)$, as a function of $t$, could be unbounded almost surely.

The remaining part of this section is devoted to the proof of the following theorem.

THEOREM 4.3. *Let $\{S_t\}_{t \geq 0}$ be the semigroup generated by (2.4) using Proposition 3.1. Assume that $c(0) > \lambda_{\mathrm{CP}}$ and $c(+\infty) < 1/\theta$. Then there exists a nontrivial probability measure $\mu$ on $(\Omega, \mathcal{B}(\Omega))$ such that $\mu[S_t f] = \mu[f]$ for all $t \geq 0$ and $f \in \mathrm{Lip}(\Omega)$.*

We need some preparatory results.

LEMMA 4.4. *Assume that $c : \mathbb{N} \to [0, +\infty)$ is a nonincreasing function such that $c(0) > \lambda_{\mathrm{CP}}$, while $c(+\infty) < 1/\theta$. For any $n \in \mathbb{N}$, define $c_n := c \mathbb{1}_{[0, n-1]}$ and consider the generator*

$$(\mathcal{L}_n f)(\eta) := \sum_x \eta(x) \left[ (\partial_x^- f)(\eta) + \sum_y c_n(\eta(y)) p(x,y) (\partial_y^+ f)(\eta) \right].$$

*Then there exists $\mu_n$ probability measure on $\Omega$ such that:*

  (i) $\mu_n \mathcal{L}_n \equiv 0$;
  (ii) *the sequence $\{\mu_n : n \geq 1\}$ is nondecreasing with respect to the stochastic ordering of measures;*
  (iii) *denote by $\nu_\lambda$ the nontrivial invariant probability measure of the CP on $X$ with parameter $\lambda := c(0) > \lambda_{\mathrm{CP}}$ (see [8], page 265). Then $\nu_\lambda \leq \mu_n$ for any $n \geq 2$;*
  (iv) *the sequence $\{\mu_n\}_{n \in \mathbb{N}}$ is tight.*

PROOF. Notice that $\mathcal{L}_n$ is of the form (2.4) so it generates a Markov process $\{\eta_{t,n}\}_{t \geq 0}$:

  (i) The process $\{\eta_{t,n}\}_{t \geq 0}$ is monotone because of Proposition 3.5 and Remark 3.11. If the initial condition is $n\mathbf{1}$, then by standard arguments (see [8], Chapter III, Theorem 2.3), $\eta_{t,n} \Rightarrow \mu_n$ as $t \to +\infty$. Furthermore, $\mu_n \mathcal{L}_n \equiv 0$.
  (ii) For any $n \geq 2$, by Proposition 3.5 and Remark 3.11, there exists a monotone coupling between $\{\eta_{t,n}\}_{t \geq 0}$, starting from $n\mathbf{1}$, and $\{\eta_{t,n+1}\}_{t \geq 0}$, starting from $(n+1)\mathbf{1}$, such that $\eta_{t,n} \leq \eta_{t,n+1}$ for any $t \geq 0$. Let $f : \Omega \to \mathbb{R}$ be a nondecreasing function, then $\mathbb{E}^{n\mathbf{1}}[f(\eta_{t,n})] \leq \mathbb{E}^{(n+1)\mathbf{1}}[f(\eta_{t,n+1})]$ for any $t \geq 0$. By taking the limit, as $t \to +\infty$, we get $\mu_n(f) \leq \mu_{n+1}(f)$.

(iii) By Proposition 3.5 we can couple $\{\eta_{t,2}\}_{t\geq 0}$, starting from $2\mathbf{1}$, and a supercritical CP $\{\xi_t\}_{t\geq 0}$, starting from $\mathbf{1}$, with parameter $\lambda = c(0)$ in such a way that $\xi_t \leq \eta_{t,2}$ for any $t \geq 0$.

(iv) Note that $c(+\infty) < 1/\theta$ implies $\bar{k} := \inf\{k \in \mathbb{N} : c(k) < 1/\theta\} < +\infty$. Take $n \geq \bar{k}$ and observe that by Proposition 3.5 and Remark 3.11 there exists a monotone coupling between $\{\eta_{t,n}\}_{t\geq 0}$, and the BRW $\{\zeta_t\}_{t\geq 0}$ generated by (3.1), with $k := \bar{k}$, $\gamma = 1$, and parameter $c(\bar{k})$, both starting from $n\mathbf{1}$. Since $\eta_{t,n} \leq \zeta_t$ for any $t \geq 0$, then $\mathbb{E}^{n\mathbf{1}}[\eta_{t,n}(0)] \leq \mathbb{E}^{n\mathbf{1}}[\zeta_t(0)]$. By taking the limit as $t \to +\infty$ and using Lemma 3.3 and Remark 3.11, we have that $\mu_n(\eta(x)) \leq U_{1,X}(\bar{k}, c(\bar{k}), 1)$ for any $n \geq \bar{k}, x \in X$. Hence, there exists a constant $C := U_{1,X}(\bar{k}, c(\bar{k}), 1)$ such that, for any $r > 0$ and $n \geq \bar{k}$, we have, by the Chebyshev inequality, $\mu_n(\eta \in \Omega : \eta(x) > r) \leq C/r$. Let us fix $A > 0$ and define $r(x) := A/\sqrt{\alpha(x)}$ for any $x \in X$. We have

$$\mu_n(\eta \in \Omega : \eta(x) \leq r(x) \text{ for any } x \in X)$$
$$= 1 - \mu_n(\eta \in \Omega : \text{there exists } x \in X \text{ such that } \eta(x) > r(x))$$
$$\geq 1 - \sum_x \mu_n(\eta \in \Omega : \eta(x) > r(x))$$
$$\geq 1 - C \sum_x \frac{1}{r(x)} = 1 - \frac{C}{A} \sum_x \sqrt{\alpha(x)}.$$

By our assumptions, $\sum_x \sqrt{\alpha(x)} < +\infty$, whence, for any $\epsilon > 0$, we can choose $A$ so that $\mu_n(\eta \in \Omega : \eta(x) \leq r(x) \text{ for any } x \in X) \geq 1 - \epsilon$ for any $n \geq \bar{k}$. The subset $K := \{\eta \in \Omega : \eta(x) \leq r(x) \text{ for any } x \in X\}$ of $\Omega$ is compact. In fact, $K = \prod_x [0, r(x)]$ since $\eta \in \prod_x ([0, r(x)] \cap \mathbb{N})$ implies that

$$\sum_x \eta(x)\alpha(x) \leq A \sum_x \sqrt{\alpha(x)} < +\infty,$$

that is, $\eta \in \Omega$. $\square$

Since the sequence $\{\mu_n\}_{n \in \mathbb{N}}$ is tight and monotone and since the set of continuous, monotone functions separates the set of probability measures, then the sequence converges weakly to a probability measure on $\Omega$, say, $\mu$.

Moreover, $\mu$ inherits all the symmetries of $X$ and $P$: if $T$ is a bijection of $X$ onto itself such that $x \sim y$ if and only if $Tx \sim Ty$ and $p(Tx, Ty) = p(x, y)$, then $\mu(T^{-1}(A)) = \mu(A)$ for all measurable set of configurations $A$. In particular, if $P$ is the simple random walk on $\mathbb{Z}^d$ or on the homogeneous tree of degree $n$, then $\mu$ is translation invariant.

By the previous lemma, $\mu_n \leq \mu$ for any $n \in \mathbb{N}$. Furthermore, $0 < \nu_\lambda(\eta(x)) \leq \mu(\eta(x)) \leq U_{1,X}(\bar{k}, c(\bar{k}), 1)$ (by the same bound on $\mu_n$), hence, $\mu$ is not $\delta_\mathbf{0}$ and $\mu[\|\eta\|] < +\infty$. We prove that $\mu$ is invariant by showing that $\mu[\mathcal{L}f] = 0$ for any $f \in \text{Lip}(\Omega)$ (see Proposition 3.13). In order to see this, we need a preparatory lemma, indeed, in the proof of Proposition 4.6 we need that $\mu_n(\mathcal{L}f) \to \mu(\mathcal{L}f)$ as $n \to +\infty$, but this does not follow directly from $\mu_n \Rightarrow \mu$ because $\mathcal{L}f$ is unbounded.

LEMMA 4.5. *Let $\{\mu_n\}_{n\in\mathbb{N}}$ be a nondecreasing sequence of probability measures on $\Omega$ and assume that $\mu_n \Rightarrow \mu$ as $n \to +\infty$ and $\mu[\|\eta\|] < +\infty$. For any $m \in \mathbb{N}$, $\eta \in \Omega$ and $g:\Omega \to \mathbb{R}$, define the configuration $\widetilde{\eta}^m(\cdot) := \mathbb{1}_{B(x_0,m)}(\cdot)\eta(\cdot)$ and $\widetilde{g}^m(\eta) := g(\widetilde{\eta}^m)$. Assume that $g:\Omega \to \mathbb{R}$ satisfies the following:*

(1) *there exists $C > 0$ such that $|g(\eta)| \le C(\|\eta\| + 1)$ for any $\eta \in \Omega$;*
(2) *$\mu_n[|g - \widetilde{g}^m|] \to 0$ as $m \to +\infty$ uniformly in $n \in \mathbb{N}$.*

*Then $\mu_n(g) \to \mu(g)$ as $n \to +\infty$.*

PROOF. We have that
$$|\mu_n[g] - \mu[g]| \le |\mu_n[g - \widetilde{g}^m]| + |\mu_n[\widetilde{g}^m] - \mu[\widetilde{g}^m]| + |\mu[\widetilde{g}^m - g]|.$$

By hypothesis (2) and the dominated convergence theorem, the first and last term on the right-hand side of the above inequality may be made small uniformly in $n$ by taking $m$ sufficiently large. Fix $m \in \mathbb{N}$ such that these terms are smaller than $\epsilon > 0$ for any $n \in \mathbb{N}$. For the middle term, define $g^{m,k}(\eta) := g^m(\eta)\mathbb{1}_{(-\infty,k]}(|g^m(\eta)|)$, $k \in \mathbb{N}$ and write

$$\begin{aligned}(4.1)\quad &|\mu_n[\widetilde{g}^m] - \mu[\widetilde{g}^m]| \\ &\qquad \le |\mu_n[\widetilde{g}^m - \widetilde{g}^{m,k}]| + |\mu_n[\widetilde{g}^{m,k}] - \mu[\widetilde{g}^{m,k}]| + |\mu[\widetilde{g}^{m,k} - \widetilde{g}^m]|.\end{aligned}$$

Note that by hypothesis (1) and elementary bounds,
$$\begin{aligned}|\widetilde{g}^m(\eta) - \widetilde{g}^{m,k}(\eta)| &= |\widetilde{g}^m(\eta)|\mathbb{1}_{(k,+\infty)}(|\widetilde{g}^m(\eta)|) \\ &\le C(\|\eta\| + 1)\mathbb{1}_{(k/C-1,+\infty)}(\|\eta\|) := v_k(\eta).\end{aligned}$$

By monotonicity, the first and the last term on the right-hand side of (4.1) can be bounded above by $\mu[v_k]$. Furthermore, $\lim_{k \to +\infty} v_k = 0$ and $v_k(\eta) \le C(\|\eta\| + 1)$, so, by dominated convergence, the first and the last term on the right-hand side of (4.1) can be made smaller than $\epsilon$ by taking $k$ large. Finally, fix $k \in \mathbb{N}$ large enough, and observe that the middle term on the right-hand side of (4.1) goes to 0 as $n \to +\infty$ by weak convergence. $\square$

PROPOSITION 4.6. *Let $\mu$ be the weak limit of the sequence $\{\mu_n\}_{n\in\mathbb{N}}$ defined in Lemma 4.4, then $\mu[\mathcal{L}f] = 0$ for any $f \in \mathrm{Lip}(\Omega)$.*

PROOF. We start splitting
$$(4.2)\qquad |\mu[\mathcal{L}f]| \le |\mu[\mathcal{L}f] - \mu_n[\mathcal{L}f]| + \mu_n[|\mathcal{L}f - \mathcal{L}_n f|].$$

Roughly speaking, the first one of these two terms goes to 0 by weak convergence, while the second one goes to 0 since $\mathcal{L}_n f \to \mathcal{L} f$.

By Lemma 4.5, $|\mu[\mathcal{L}f] - \mu_n[\mathcal{L}f]| \to 0$ if we can show that condition (2) is satisfied by $g := \mathcal{L}f$ [condition (1) is easily verified]. Observe that

$$
\begin{aligned}
|(\mathcal{L}f)(\eta) - (\mathcal{L}f)(\widetilde{\eta}^m)| \\
\leq \sum_x \mathbb{1}_{B(x_0,m)}(x)\eta(x) & \left[|(\partial_x^- f)(\eta) - (\partial_x^- f)(\widetilde{\eta}^m)|\right. \\
& \left. + c(0) \sum_y p(x,y)|(\partial_y^+ f)(\eta) - (\partial_y^+ f)(\widetilde{\eta}^m)|\right] \\
+ \sum_x \mathbb{1}_{B(x_0,m)^{\complement}}(x)\eta(x) & \left[|(\partial_x^- f)(\eta)| + c(0) \sum_y p(x,y)|(\partial_y^+ f)(\eta)|\right].
\end{aligned}
\tag{4.3}
$$

It easy to show that $L(\partial_z^\pm f) \leq 2L(f)$ for any $f \in \text{Lip}(\Omega)$, $x \in X$. Thus, the first of the two sums on the right-hand side of (4.3) is dominated by

$$2L(f)(1+c(0))\|\eta - \widetilde{\eta}^m\| \sum_x \mathbb{1}_{B(x_0,m)}(x)\eta(x).$$

Moreover,

$$\|\eta - \widetilde{\eta}^m\| \sum_x \mathbb{1}_{B(x_0,m)}(x)\eta(x) = \sum_{x,y} \alpha(y)\mathbb{1}_{B(x_0,m)}(x)\mathbb{1}_{B(x_0,m)^{\complement}}(y)\eta(x)\eta(y).$$

This implies that the $\mu_n$ mean of the first term on the right-hand side of (4.3) is smaller or equal to

$$2L(f)(1+c(0)) \sum_{x,y} \mathbb{1}_{B(x_0,m)}(x)\mathbb{1}_{B(x_0,m)^{\complement}}(y)\alpha(y)\mu_n[\eta(x)\eta(y)].$$

This term goes to 0 uniformly in $n \in \mathbb{N}$ since $\mu_n[\eta(x)\eta(y)] \leq \mu[\eta(x)\eta(y)] \leq C$ by Remark 3.11 and Lemma 3.3 while, by our choice of $M > (D-1)^2$,

$$\sum_{x,y} \mathbb{1}_{B(x_0,m)}(x)\mathbb{1}_{B(x_0,m)^{\complement}}(y)\alpha(y) \stackrel{m \to \infty}{\longrightarrow} 0.$$

The second term on the right-hand side of (4.3) can be dominated by

$$L(f)(1+c(0)M) \sum_x \mathbb{1}_{B(x_0,m)^{\complement}}(x)\eta(x)\alpha(x),$$

hence, its $\mu_n$ mean converges uniformly in $n$ to 0, since it is not larger than

$$U_{1,X}(\bar{k}, c(\bar{k}), \theta) \sum_x \mathbb{1}_{B(x_0,m)^{\complement}}(x)\alpha(x) \stackrel{m \to \infty}{\longrightarrow} 0$$

(again use Lemma 3.3).

We are left with the proof that the second term on the right-hand side of (4.2) goes to 0 as $n \to +\infty$. Observe that since $\mu_n$ is concentrated on $\{\eta : \eta \leq n\mathbf{1}\}$,

$$\mu_n[|\mathcal{L}f - \mathcal{L}_n f|] = c(n) \sum_{x,y} p(x,y) \mu_n[\eta(x) \mathbb{1}_{\{n\}}(\eta(y)) |\partial_y^+ f(\eta)|]$$

$$\leq c(0) L(f) \sum_{x,y} p(x,y) \alpha(y) \mu_n[\eta(x) \mathbb{1}[n, +\infty)(\eta(y))].$$

By the Schwarz and Chebyshev inequalities, Lemma 3.3 and Remark 3.11, we have that

$$\mu_n[\eta(x) \mathbb{1}_{[n,+\infty)}(\eta(y))] \leq \sqrt{\mu_n[\eta(x)^2] \mu_n(\eta(y) \geq n)}$$

$$\leq \sqrt{\mu[\eta(x)^2] \mu(\eta(y) \geq n)} \leq \frac{C}{\sqrt{n}} \to 0,$$

where $C$ does not depend on $x$ and $y$. □

PROOF OF THEOREM 4.3. It follows from Lemma 4.4, Lemma 4.5 and Proposition 4.6. □

**5. Final remarks and examples.** In this paper we consider mainly *local survival* (i.e., persistence of the population in a fixed site). We already observed that, for general $P$ (e.g., for the simple random walk on a general graph), a weaker type of survival (the weak phase) is possible for both the CP and the BRW. One can associate to this *global phase* two critical parameters (which coincide with $\lambda_{\text{CP}}$ and $\lambda_{\text{BRW}}$, for instance, on $\mathbb{Z}^d$). Clearly this phenomenon could be observed also in the evolution of a RBRW.

Proposition 4.1 does not describe the behavior of all possible RBRW: in particular, it is not clear what happens if $1/\rho < c(0) \leq \lambda_{\text{CP}}$ or if $c(+\infty) \in (1/\theta, 1/\rho)$. As for this last question, one would ask whether this interval may be nonempty. We already noticed that $\rho = \theta$ on amenable graphs with a strongly reversible random walk (such as the simple random walk on $\mathbb{Z}^d$), nevertheless, there are examples of graphs (see Example 5.1) where $\rho \neq \theta$ and for any $\lambda \in (1/\theta, 1/\rho)$, the BRW starting from a finite configuration vanishes locally with probability 1, while starting from a bounded configuration greater than $\mathbf{1}$, the expected number of individuals at a fixed site diverges. Roughly speaking, this is possible on graphs where the contribution of far distant individuals is not negligible; indeed, if $P$ is symmetric, this behavior is equivalent to the existence of a weak phase for the BRW. One may conjecture that the critical parameter of this phenomenon may be $1/\limsup_{n \to \infty} \sqrt[n]{\sum_x p^{(n)}(x,y)}$ (which does not depend on $y$). Finally, another open question is the extremality of the invariant measure $\mu$.

EXAMPLE 5.1. Let us consider the homogeneous tree $\mathbb{T}_{n+1}$ where the degree of each vertex is $n+1$ and choose $n \geq 2$. Fix a reference vertex $o$ and $p \in [0, 1/n]$. Given two neighbors $x$ and $y$, we define $p(x,y)$ as $p$ if $|x|+1 = |y| \geq 2$ (recall that $|x|$ is the distance from $x$ to $o$), as $1/(n+1)$ if $y \sim x = o$ and as $1 - np$ otherwise. By using standard generating function techniques and the fact that $\|P^T\| \geq \theta \geq \limsup_{n \to \infty} \sqrt[n]{\sum_y p^{(n)}(y,x)} \geq \rho$ (for all $x$), it is easy to show that Note that $\rho < 1$ for all $p > 1/2n$ and that $\rho < \theta$,

| Range of $p$ | $\rho$ | $\theta$ |
|---|---|---|
| $[0, 1/2n]$ | $1$ | $1 < n - (n^2-1)p \leq \theta \leq (n+1)(1-np)$ |
| $(1/2n, 1/(n+1))$ | $2\sqrt{np(1-np)}$ | $1 < n - (n^2-1)p \leq \theta \leq (n+1)(1-np)$ |
| $1/(n+1)$ | $2\sqrt{n}/(n+1)$ | $\theta = 1$ |
| $[1/(n+1), 1/n]$ | $2\sqrt{np(1-np)}$ | $\theta = n - (n^2-1)p < 1$ |

for instance, if $p = 1/(n+1)$, that is, the simple random walk. In this case $\theta = 1$; by using the explicit solution given in Remark 3.4 (with $\varphi \equiv 1$), and noting that, by translation invariance, $\sum_x p^{(n)}(x,y) = \|(P^T)^n\| \geq \theta^n = 1$, we have that $\lim_{t \to \infty} \mathbb{E}^\varphi[\eta_t(x)] = +\infty$ if $\lambda > 1$. Moreover, considering the BRW with at least one immortal particle per site, if $\lambda > 1$, then the expected number of individuals on a fixed site diverges as $t$ goes to infinity.

**Acknowledgment.** The authors wish to thank Paolo Dai Pra for his invaluable suggestions.

D. BERTACCHI
DIPARTIMENTO DI MATEMATICA E APPLICAZIONI
UNIVERSITÀ DI MILANO–BICOCCA
VIA COZZI 53
20125 MILANO
ITALY
E-MAIL: daniela.bertacchi@unimib.it

G. POSTA
F. ZUCCA
DIPARTIMENTO DI MATEMATICA
POLITECNICO DI MILANO
PIAZZA LEONARDO DA VINCI 32
20133 MILANO
ITALY
E-MAIL: gustavo.posta@polimi.it
    fabio.zucca@polimi.it